\documentclass[10pt]{article}
\usepackage{times}

\DeclareMathVersion{varnormal}
\DeclareMathVersion{varbold}

\SetSymbolFont{letters}{normal}{OML}{txmi}{m}{it}
\SetSymbolFont{letters}{bold}{OML}{txmi}{bx}{it}
\SetSymbolFont{letters}{varnormal}{OML}{mdugm}{m}{it}
\SetSymbolFont{letters}{varbold}{OML}{mdugm}{b}{it}
\SetSymbolFont{operators}{normal}{OT1}{txr}{m}{n}
\SetSymbolFont{operators}{bold}{OT1}{txr}{bx}{n}
\SetSymbolFont{operators}{varnormal}{OT1}{mdugm}{m}{n}
\SetSymbolFont{operators}{varbold}{OT1}{mdugm}{b}{n}
\SetSymbolFont{symbols}{normal}{OMS}{txsy}{m}{n}
\SetSymbolFont{symbols}{bold}{OMS}{txsy}{bx}{n}
\SetSymbolFont{symbols}{varnormal}{OMS}{mdugm}{m}{n}
\SetSymbolFont{symbols}{varbold}{OMS}{mdugm}{b}{n}
\SetSymbolFont{largesymbols}{normal}{OMX}{txex}{m}{n}
\SetSymbolFont{largesymbols}{bold}{OMX}{txex}{bx}{n}
\SetSymbolFont{largesymbols}{varnormal}{OMX}{mdugm}{m}{n}
\SetSymbolFont{largesymbols}{varbold}{OMX}{mdugm}{b}{n}

\SetMathAlphabet{\mathrm}{normal}{OT1}{txr}{m}{n}
\SetMathAlphabet{\mathrm}{bold}{OT1}{txr}{bx}{n}
\SetMathAlphabet{\mathrm}{varnormal}{OT1}{mdugm}{m}{n}
\SetMathAlphabet{\mathrm}{varbold}{OT1}{mdugm}{b}{n}

\SetMathAlphabet{\mathit}{normal}{OT1}{txr}{m}{it}
\SetMathAlphabet{\mathit}{bold}{OT1}{txr}{bx}{it}
\SetMathAlphabet{\mathit}{varnormal}{OT1}{mdugm}{m}{it}
\SetMathAlphabet{\mathit}{varbold}{OT1}{mdugm}{b}{it}

\usepackage{amsmath}
\usepackage{amsthm}\theoremstyle{plain}
\usepackage{mathrsfs}
\usepackage{amssymb}
\usepackage[ruled]{algorithm2e}
\usepackage{program}
\usepackage{empheq}
\usepackage{subcaption}
\usepackage[paper=a4paper,top=36pt,left=100pt,right=100pt,bottom=50pt,foot=25pt]{geometry}
\usepackage{titlesec}
\usepackage{lineno}
\usepackage{hyperref}
\usepackage{nicefrac}
\usepackage{float}
\usepackage{upgreek}
\usepackage[dvipsnames]{xcolor}
\usepackage[numbers,compress]{natbib}
\usepackage{verbatim}
\usepackage{listings}
\usepackage{mdframed}
\usepackage[T1]{fontenc}

\titleformat{\section}{\normalfont\bfseries}{\thesection}{1em}{}
\titleformat{\subsection}{\normalfont}{\thesubsection}{1em}{}
\titleformat{\subsubsection}{\normalfont\it}{\thesubsubsection}{1em}{}

\renewcommand\appendix{\par
  \setcounter{section}{0}
  \setcounter{subsection}{0}
  \setcounter{figure}{0}
  \setcounter{table}{0}
  \setcounter{equation}{0}
  \renewcommand\thesection{Appendix \Alph{section}}
  \renewcommand\thefigure{\Alph{section}\arabic{figure}}
  \renewcommand\thetable{\Alph{section}\arabic{table}}
  \renewcommand\theequation{\Alph{section}\arabic{equation}}  
}

\definecolor{cmtcolor}{RGB}{128,128,128}
\definecolor{kwdcolor}{RGB}{180,20,128}
\definecolor{strcolor}{RGB}{0,120,80}

\newcommand\digitstyle{\color{kwdcolor2}}
\makeatletter
\newcommand{\ProcessDigit}[1]
{%
  \ifnum\lst@mode=\lst@Pmode\relax%
   {\digitstyle #1}%
  \else
    #1%
  \fi
}

\definecolor{kwdcolor2}{RGB}{20,70,180}
\definecolor{hltcolor}{RGB}{235,120,60}

\def\ttt#1;{\texttt{#1}} 

\newcommand{\folderPath}{\string"/dev/null\string"}

\title{\large Half precision wave simulation}

\author{
\small Longfei Gao\thanks{Email address: longfei.gao@anl.gov} \, Kevin Harms\thanks{Email address: harms@alcf.anl.gov}\\
{\it \small Argonne National Laboratory, 9700 S Cass Ave, Lemont, IL 60439} 
}
\date{}

\def\bgnEqn{\begin{equation}}
\def\endEqn{\end{equation}}

\makeatletter\newcommand{\eqrefNolink}[1]{\textup{\tagform@{\ref*{#1}}}}\makeatother

\begin{document}


\maketitle

\begin{abstract}
In recent years, half precision floating-point arithmetic has gained wide support in hardware and software stack thanks to the advance of artificial intelligence and machine learning applications. 
Operating at half precision can significantly reduce the memory footprint comparing to operating at single or double precision. 
For memory bound applications such as time domain wave simulations, this is an attractive feature. 
However, the narrower width of the half precision data format can lead to degradation of the solution quality due to larger roundoff errors.
In this work, we illustrate with carefully designed numerical experiments the negative impact caused by the accumulation of roundoff errors in wave simulations.
Specifically, the energy-conserving property of the wave equations is employed as a convenient diagnosis tool.
The corresponding remedy in the form of compensated sum is then provided, with its efficacy demonstrated using numerical examples with both acoustic and elastic wave equations on hardware that support half precision arithmetic natively.
\end{abstract}

\section{Introduction}
In traditional science and engineering applications, we are accustomed to double or single precision floating-point arithmetic thanks to the wide adoption of the IEEE standard for Floating-Point Arithmetic \cite{IEEE2019Standard} first established in the 80s. 
In recent years, because of the advance in artificial intelligence and machine learning applications, half precision floating-point arithmetic has gained popularity as well.
Modern hardware, GPUs in particular, often have significant chip area dedicated to low precision units. 
How to effectively utilize these low precision units in traditional science and engineering applications remains an open question.

For seismic wave simulations, it has long been recognized that single precision is often adequate to achieve satisfactory simulation quality; see \cite[p.14]{manual2021specfem3d} and \cite{abdelkhalek2009fast} for example. 
%
%
Additional justifications 
include noise in the field-collected data and approximations in the underlying physical model.
The benefit of operating at single precision versus double precision resides in the reduced memory footprint. 
On modern hardware, memory speed is often the limiting factor for execution time; see Figure 2.2 of \cite[p.80]{Hennessy2015computer} for a historical trend of the processor-DRAM performance gap.
By reducing the operating precision from double to single, the storage and movement requirement per datum is halved.
For memory bound applications such as seismic processing, 
we can expect the execution time to be significantly reduced, as well as substantial saving in power consumption.

However, little effort has emerged in operating at half precision for seismic processing tasks.
A notable effort can be found in \cite{fabien2020seismic}, where the author explored the performance benefit and solution quality of operating at half precision for modeling, imaging, and inversion tasks.
The results shown there is encouraging despite the accuracy issue related with the narrow data format in half precision. 

In this work, we focus on the particular issue of degradation in solution quality due to the accumulation of roundoff errors in half precision wave simulations.  
Specifically, we demonstrate using numerical examples the potential pitfalls of naively switching to half precision. 
The energy-conserving property of the wave equation is employed as a convenient diagnosis tool.
Remedy in the form of compensated sum is then provided to address these pitfalls. 
Its efficacy is demonstrated using both acoustic and elastic wave equations, with the former tested on Intel Ponte Vecchio GPU and the later on NVIDIA A100 GPU.
Additionally, a discussion section is included to bring to the readers' attention other potential issues that may occur when operating at half precision but not addressed here.


\section{Terminology}
%
The IEEE standard specifies floating point formats that are most widely supported on modern hardware, which will be the focus of this work.
The 64-bit double and 32-fit single formats are historically the most prevalent.
The 16-bit half format has also gained popularity in recent years. 
These formats are referred to as fp64, fp32, and fp16 in the following, respectively. 
Some of their relevant characteristics are summarized in Table \ref{tbl_formats}. 
%
%

The IEEE standard stipulates that addition, subtraction, multiplication, division, square root, and fused multiply add operations shall be performed as if an intermediate result correct to infinite precision and with unbounded range has been produced first, and then rounded to fit the destination format.
The addition operation will be the focus of this work.
%
With an implementation compliant to the IEEE standard, these operations will incur a relative error no larger than the unit roundoff for normal numbers.

\begin{table}[H]
\captionsetup{width=.9\linewidth,font=small}
\caption{Summary of the fp64, fp32, and fp16 formats. (In the row regarding number of fraction bits, $+1$ means the hidden bit, which is a benefit of working with the binary number system.)}
\label{tbl_formats}
\centering
\begin{tabular}{cccc}
						  \hline
                          & fp64 & fp32 & fp16 \\ \hline
number of exponent bits   & 11   & 8      & 5  \\
number of fraction bits   & 52+1 & 23+1 & 10+1 \\
unit roundoff 			  
						  & $1.1102\times 10^{-16}$
						  & $5.9605\times 10^{-8} $
						  & $4.8828\times 10^{-4} $ 
						  \\
largest number            
                          & $1.7977\times 10^{308}$ 
                          & $3.4028\times 10^{38} $ 
                          & $6.5504\times 10^{4}  $ 
                          \\
smallest positive normal number
                          & $2.2251\times 10^{-308}$
                          & $1.1755\times 10^{-38 }$
                          & $6.1035\times 10^{-5}$
                          \\
smallest  positive subnormal number
                          & $4.9407\times 10^{-324}$
                          & $1.4013\times 10^{-45 }$
                          & $5.9605\times 10^{-8}$
						  \\ \hline
\end{tabular}
\end{table}

\section{Problem description}\label{sec_problem}
To start, we consider the following 2D acoustic wave equation 
\begin{equation}
\arraycolsep=2.5pt\def\arraystretch{2.25}
\label{acoustic_wave_equation}
\left\{
\begin{array}{lcl}
\displaystyle \rho  \frac{\partial v_x}{\partial t} & = & \displaystyle \frac{\partial p}{\partial x} \, ;
\\
\displaystyle \rho  \frac{\partial v_y}{\partial t} & = & \displaystyle \frac{\partial p}{\partial y} \, ;
\\
\displaystyle \beta \frac{\partial p  }{\partial t} & = & \displaystyle \frac{\partial v_x}{\partial x} \ + \ \frac{\partial v_y}{\partial y} \ + \ s_p \, ,
\end{array}
\right.
\end{equation}
where $v_x$, $v_y$, and $p$ are the sought solution variables, standing for particle velocities in $x$ and $y$ directions and pressure, respectively; $\rho$ and $\beta$ are given physical parameters, standing for density and compressibility, respectively; $s_p$ is a given source term that excites the wave propagation.
The dependency of $v_x$, $v_y$, $p$, and $s_p$ on time $t$ is implicit in the above equation for brevity.
In this work, we consider a point source on pressure, i.e., $s_p$, which concentrates at a singular point of the domain, which has a temporal profile of Ricker wavelet with 5~Hz central frequency.
The system is assumed to be at rest for the initial state.
Periodic boundary condition is considered for both directions in the main body of this work.
Results from the elastic wave equation with free surface boundary condition can be found in \ref{App_elastic}.

\begin{figure}[H]
\centering
\includegraphics[scale=0.05]{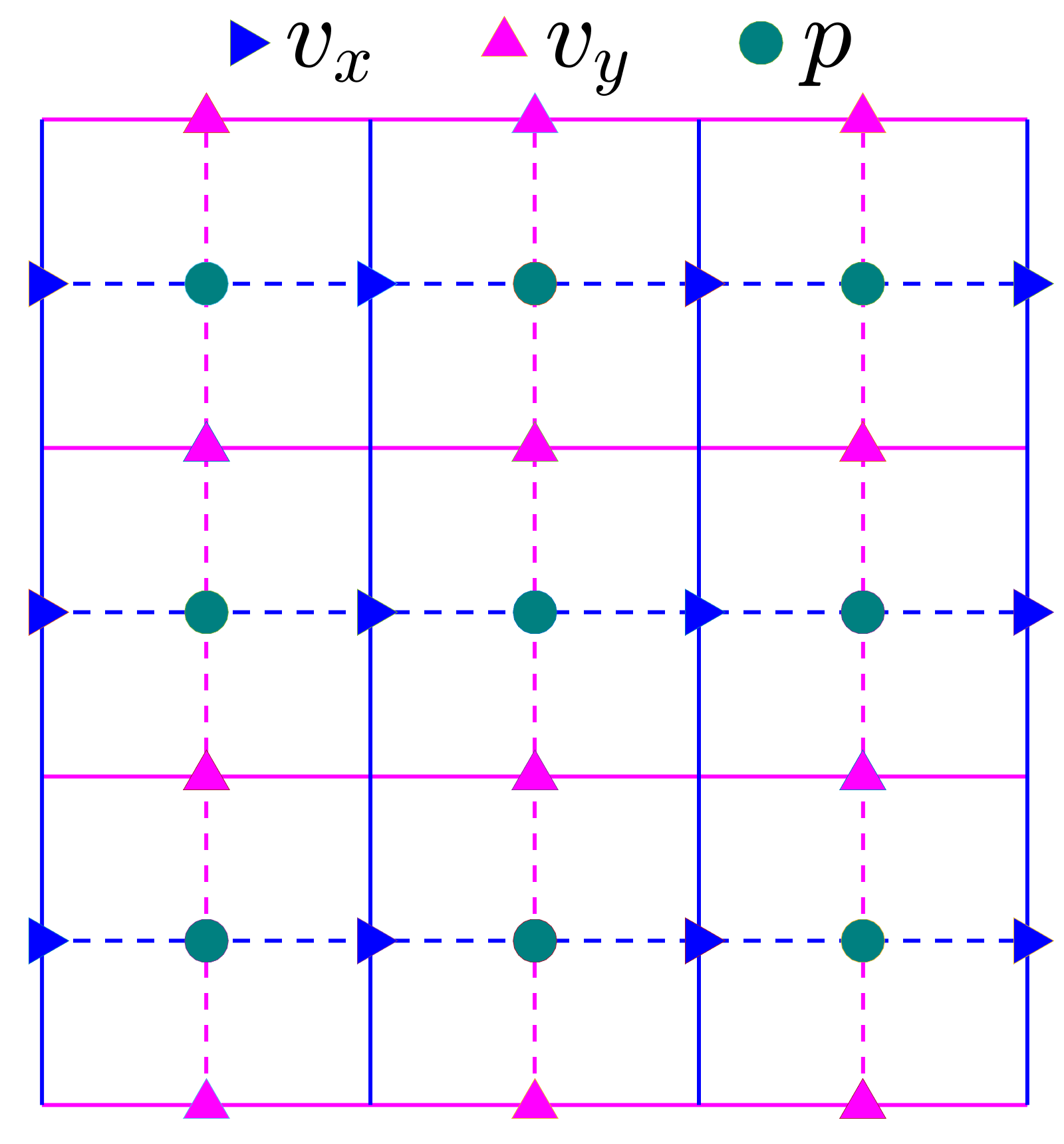}
\caption{Illustration of the grid layout in space.}
\label{space_layout}
\end{figure}

We discretize equation \eqref{acoustic_wave_equation} in space on a set of staggered grids as illustrated in Figure \ref{space_layout}. 
The fourth order finite difference stencil 
$
[\nicefrac{1}{24}, -\nicefrac{9}{8}, \nicefrac{9}{8}, -\nicefrac{1}{24}]
$
(see \cite{levander1988fourth})
is used to approximate the spatial derivatives $\frac{\partial}{\partial x}$ and $\frac{\partial}{\partial y}$.
Periodic boundary condition is imposed by wrapping the stencil around the domain. 

\begin{figure}[H]
\centering
\includegraphics[scale=0.05]{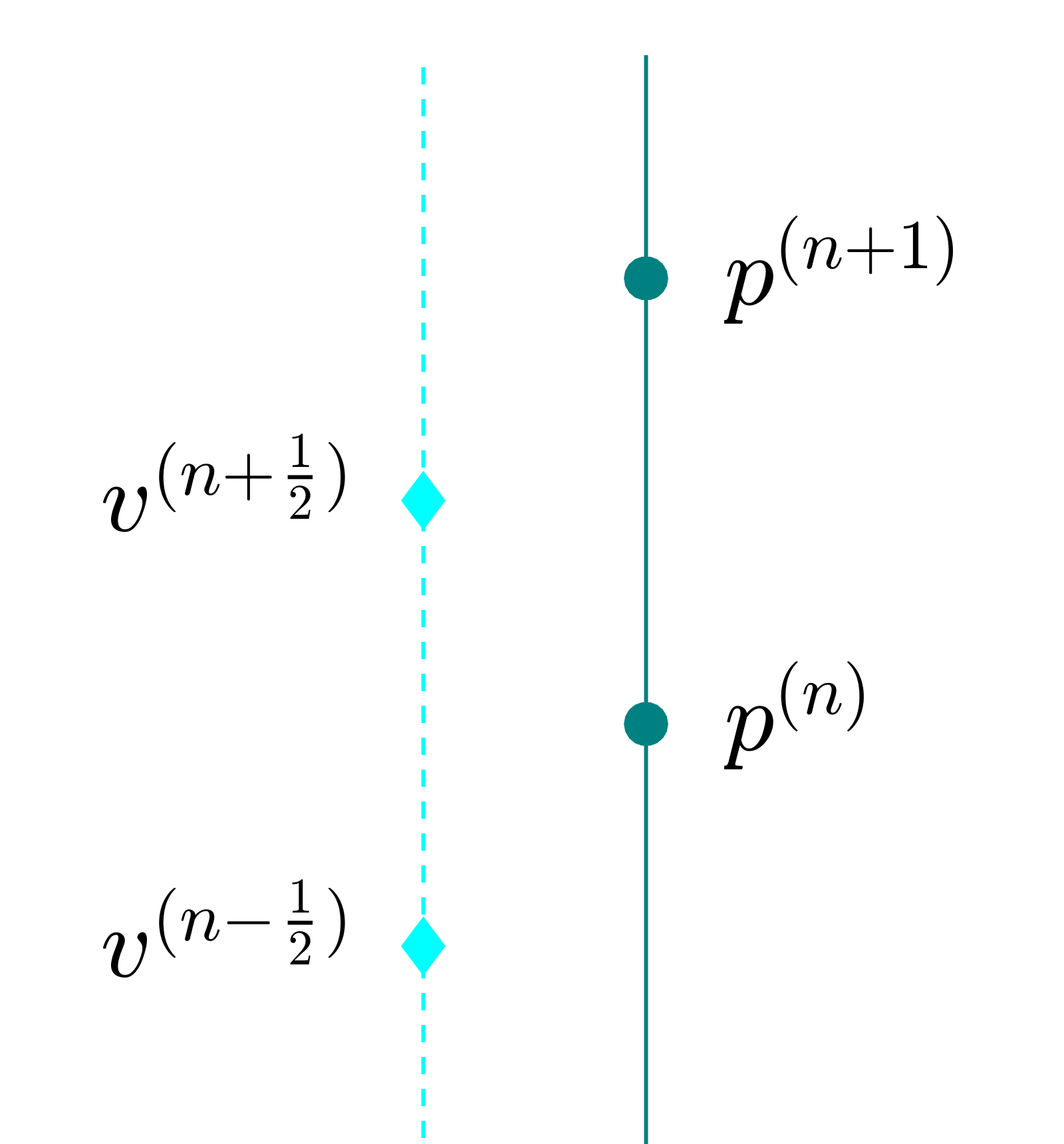}
\caption{Illustration of the grid layout in time.}
\label{time_layout}
\end{figure}

Temporal discretization of equation \eqref{acoustic_wave_equation} is achieved via the second order staggered leapfrog scheme. 
Layout of the discretized variables in time is illustrated in Figure \ref{time_layout}.

With the above choices on discretization, the discrete wave system preserves the energy conserving property, which is a useful diagnosis tool for our experiments and precludes instability that stems from the kind of discretization error illustrated in \cite[Figure 6]{gao2018long} so that we can focus on the effect of roundoff errors accumulated during the simulation.
For more detail on the energy definition and property of the discrete wave system, the readers are referred to \cite{gao2020explicit}.

\begin{figure}[H]
\centering
\includegraphics[scale=0.225]{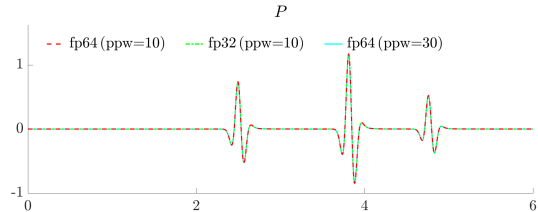}
\captionsetup{width=.9\linewidth,font=small}
\caption{Time history of $P$ at the receiver location.}
\label{fp64_fp32_P}
\end{figure}

\begin{figure}[H]
\centering
\includegraphics[scale=0.225]{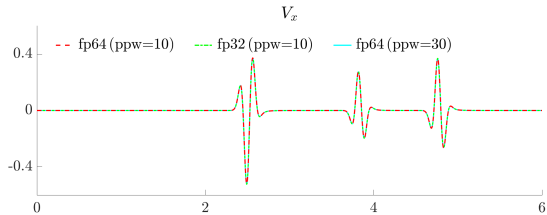}
\captionsetup{width=.9\linewidth,font=small}
\caption{Time history of $V_x$ at the receiver location.}
\label{fp64_fp32_Vx}
\end{figure}

\begin{figure}[H]
\centering
\includegraphics[scale=0.225]{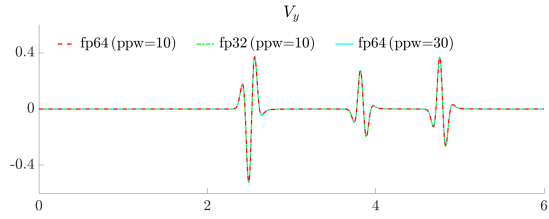}
\captionsetup{width=.9\linewidth,font=small}
\caption{Time history of $V_y$ at the receiver location.}
\label{fp64_fp32_Vy}
\end{figure}

\begin{figure}[H]
\centering
\includegraphics[scale=0.225]{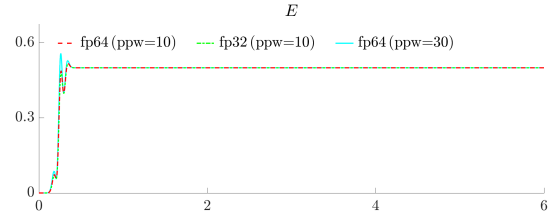}
\captionsetup{width=.9\linewidth,font=small}
\caption{Time history of the overall energy.}
\label{fp64_fp32_E}
\end{figure}

Using the discretization described above, we obtain fp64 and fp32 simulation results illustrated in Figures \ref{fp64_fp32_P} - \ref{fp64_fp32_E} for the simulation configuration outlined in the following.
The baseline spatial grid consists of 600 grid spacings on both $x$ and $y$ directions, with the minimal wavelength spanning 10 grid spacings (denoted as ppw=10 in the following). 
Homogeneous medium with $\rho=1$~$\text{kg}/\text{m}^3$ and $c=1$~$\text{m}/\text{s}$ is considered.
The maximal frequency in the source content is counted as 12.5 Hz, which leads to $\Delta x = \frac{1}{125}$~m for the baseline grid with the above specification.
Time step length $\Delta t$ is chosen as 1e-4~s.
Number of time steps is chosen as 60000.
The point source is placed at one-third of the domain on both directions.
A receiver is placed at two-thirds of the domain on both directions.

The experiments in this and the next section are tested on Intel Ponte Vecchio GPU, which supports fp16 arithmetic natively, and implemented using SYCL, where the fp16 data type can be accessed through the type \texttt{sycl::half} or the type \texttt{\_Float16} from the intel \texttt{icpx} compiler.

\begin{figure}[H]
\centering
\includegraphics[scale=0.225]{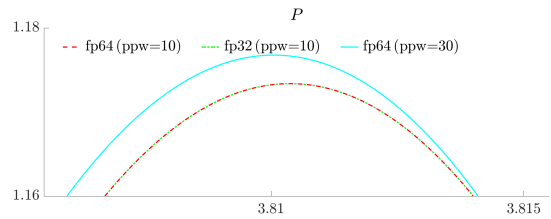}
\captionsetup{width=.9\linewidth,font=small}
\caption{Zoom-in plot of Figure \ref{fp64_fp32_P} reveals that the two simulations fp64 (ppw=10) and fp32 (ppw=10) are closer to each other than to the fp64 (ppw=30) simulation, indicating that the dominant source of error in the fp32 (ppw=10) simulation is from discretization rather than from rounding.}
\label{fp64_fp32_P_zoom_in}
\end{figure}

Time histories of the solution variables recorded at the receiver are shown in Figures \ref{fp64_fp32_P} - \ref{fp64_fp32_Vy}. 
Time history of the overall energy is shown in Figure \ref{fp64_fp32_E}.
Results from three different simulations, fp64 (ppw=10), fp32 (ppw=10), and fp64 (ppw=30) are shown in each figure.
All three simulations agree with one another nicely.
The agreement between fp64 (ppw=10) and fp64 (ppw=30) indicates that our code implementation is reliable for the ensuing experiments.
The agreement between fp64 (ppw=10) and fp32 (ppw=10) confirms the observation from earlier works that fp32 is typically sufficient for this application.
In fact, zoom-in plot of the recorded pressure, shown in Figure \ref{fp64_fp32_P_zoom_in}, reveals that the simulation results from fp64 (ppw=10) and fp32 (ppw=10) are closer to each other than to fp64 (ppw=30), indicating that discretization error dominates roundoff error in the case of fp32 (ppw=10) simulation.
Finally, from Figure \ref{fp64_fp32_E} we observe that the energy remains flat after the source effect tapers off in these simulations, reflecting the energy-conserving property in the wave system.

\renewcommand{\folderPath}{\string"\string./figure/fp64_and_fp32/\string"}

\begin{figure}[H]
\centering
\includegraphics[scale=0.225]{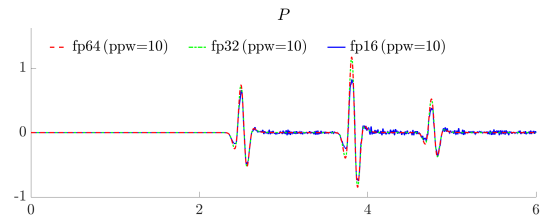}
\captionsetup{width=.9\linewidth,font=small}
\caption{Time history of $P$ at the receiver location.}
\label{fp16_naive_P}
\end{figure}

\begin{figure}[H]
\centering
\includegraphics[scale=0.225]{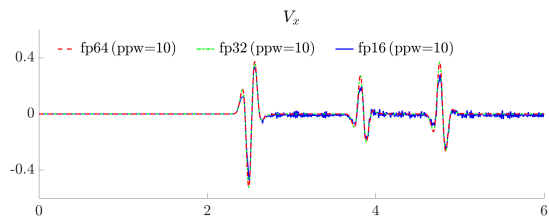}
\captionsetup{width=.9\linewidth,font=small}
\caption{Time history of $V_x$ at the receiver location.}
\label{fp16_naive_Vx}
\end{figure}

\begin{figure}[H]
\centering
\includegraphics[scale=0.225]{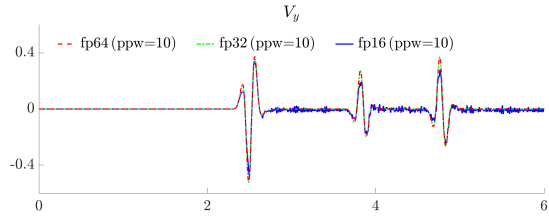}
\captionsetup{width=.9\linewidth,font=small}
\caption{Time history of $V_y$ at the receiver location.}
\label{fp16_naive_Vy}
\end{figure}

\begin{figure}[H]
\centering
\includegraphics[scale=0.225]{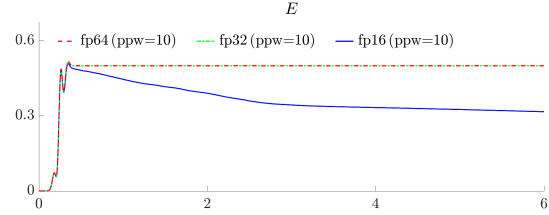}
\captionsetup{width=.9\linewidth,font=small}
\caption{Time history of the overall energy.}
\label{fp16_naive_E}
\end{figure}

Next, we reduce the operating precision to fp16 and show the simulation results in Figures \ref{fp16_naive_P} - \ref{fp16_naive_E} with comparisons to the fp64 and fp32 results.
From these figures, we observe that the fp16 simulation quality is degraded, with undesirable wiggles (Figures \ref{fp16_naive_P} - \ref{fp16_naive_Vy}) and energy loss (Figure \ref{fp16_naive_E}).
This is not surprising given the large unit roundoff ($4.8828 \times 10^{-4}$) associated with fp16.
In the next section, we outline where these issues originate from and present the remedy accordingly.

\renewcommand{\folderPath}{\string"\string./figure/fp16_naive/\string"}

\section{Compensated sum} 
\label{sec_compensated_sum}

Algorithm \ref{alg_basic} describes the procedure used in the above simulations, where lines 1, 4, 7, and 8 represent the stencil applications associated with $\frac{\partial}{\partial x}$ and $\frac{\partial}{\partial y}$ from equation \ref{acoustic_wave_equation}; lines 3, 6, and 11 stem from the solution updates associated with $\frac{\partial}{\partial t}$ from equation \ref{acoustic_wave_equation}.
%
%
We have identified that the addition operations at solution updates, 
which form a recursive sum in disguise, 
are the major contributor to the issues in half precision simulations illustrated in the previous section. 
Interest readers may refer to \ref{App_promotion} for experiments illustrating that the stencil applications are not a significant contributor to those issues, despite involving majority of the floating point operations.

One technique to improve the accuracy of sum is the so called compensated sum, which dates back to the early years of computing when double precision arithmetic was not universally available. 
The general idea is to keep track of the lost bits in floating point addition that are due to finite precision and apply (compensate) it to the sum later. 
Two common versions of such compensated sum algorithms are presented in Code snippets \ref{sum_3op} and \ref{sum_6op}, which take 3 and 6 operations, respectively, and are referred to as \texttt{sum\_3op} and \texttt{sum\_6op} in the following. 
%
%
These algorithms are often credited to \cite{kahan1965pracniques,moller1965quasi}.
%

\renewcommand{\lstlistingname}{Code}
\renewcommand{\lstlistlistingname}{List of \lstlistingname s}

\begin{center}
\begin{minipage}{0.5\linewidth}
\captionsetup[lstlisting]{font={small}}
\begin{lstlisting}[tabsize=4,basicstyle=\ttfamily\footnotesize,xleftmargin=1.em,caption={Compensated sum (3op version).},label={sum_3op}]
template<typename T>
void sum_3op ( T const a , T const b ,
               T &     s , T &     t )
{
      s = a + b;
    T z = s - a;
      t = b - z;
}
\end{lstlisting}
\end{minipage}
\end{center}

\begin{center}
\begin{minipage}{0.5\linewidth}
\captionsetup[lstlisting]{font={small}}
\begin{lstlisting}[tabsize=4,basicstyle=\ttfamily\footnotesize,xleftmargin=1.em,caption={Compensated sum (6op version).},label={sum_6op}]
template<typename T>
void sum_6op ( T const a , T const b ,
               T &     s , T &     t )
{
	    s =   a + b;
	T p_a =   s - b;
	T p_b =   s - p_a;

	T d_a =   a - p_a;
	T d_b =   b - p_b;

	    t = d_a + d_b;
}
\end{lstlisting}
\end{minipage}
\end{center}

Both \texttt{sum\_3op} and \texttt{sum\_6op} take two input arguments \texttt{a} and \texttt{b} and return two output arguments \texttt{s} and \texttt{t}, with \texttt{s} being the floating point sum of \texttt{a} and \texttt{b} and \texttt{t} keeping track of the lost bits when producing \texttt{s}.
These algorithms are sometimes referred to as error free transformations in the sense that \texttt{a} + \texttt{b} = \texttt{s} + \texttt{t} under suitable conditions (`+' here means exact sum).

The difference between \texttt{sum\_3op} and \texttt{sum\_6op} is that \texttt{sum\_3op} works only when the exponent of \texttt{a} is no less than the exponent of \texttt{b}. (For binary numbers, this is equivalent to \string|\texttt{a}\string|$\geq$\string|\texttt{b}\string|.)
Although adding an \texttt{if} statement followed by conditionally swapping the two input arguments can guard against the case \string|\texttt{a}\string|$<$\string|\texttt{b}\string|, branching is highly undesirable on modern hardware.
Instead, \texttt{sum\_6op} works for both cases by applying additional floating point operations, which is more suitable on modern hardware since the cost of these additional floating point operations is often negligible compared to, e.g., memory operations.
Nevertheless, we observe little difference in the simulation quality when applying \texttt{sum\_3op} or \texttt{sum\_6op}. This is likely because although \texttt{sum\_3op} is not an error free transformation when \string|\texttt{a}\string|$<$\string|\texttt{b}\string|, the error is still well bounded when it is applied to the sum of a sequence of numbers; see \cite[Theorem 8]{goldberg1991every} and \cite[P.~229 and p.~572]{knuth1981art}.




\begin{minipage}{0.425\linewidth}
\vspace{1em}
\begin{algorithm}[H]
\footnotesize
\SetAlCapHSkip{0ex}
\caption{Baseline simulation procedure.} 
\label{alg_basic}
\For {$i_t = 0 \cdots N_t$}
{
	\For {$i_s = 0 \cdots N^{v_x}_s$}
	{
		$
		\arraycolsep=1pt\def\arraystretch{1.75}
		\begin{array}{rrl}
\nl\hfill	{R^{v_x}}{[i_s]} & =& \frac{\partial p}{\partial x}{[i_s]} \\
\nl\hfill	{R^{v_x}}{[i_s]} &/\!=& {\rho^{v_x}}{[i_s]} \\
\nl\hfill	   {V_x}{[i_s]} &+\!\!=& {R^{v_x}}{[i_s]}
		\end{array}
		$
	} 

	\For {$i_s = 0 \cdots N^{v_y}_s$}
	{
		$
		\arraycolsep=1pt\def\arraystretch{1.75}
		\begin{array}{rrl}
\nl\hfill	{R^{v_y}}{[i_s]} & =& \frac{\partial p}{\partial y}{[i_s]} \\
\nl\hfill	{R^{v_y}}{[i_s]} &/\!=& {\rho^{v_y}}{[i_s]} \\
\nl\hfill	   {V_y}{[i_s]} &+\!\!=& {R^{v_y}}{[i_s]}
		\end{array}
		$
	}

	\For {$i_s = 0 \cdots N^{p}_s$}
	{
		$
		\arraycolsep=1pt\def\arraystretch{1.75}
		\begin{array}{rrl}
\nl\hfill	{R^{p}}{[i_s]} & =& \frac{\partial v_x}{\partial x}{[i_s]} \\
\nl\hfill	{R^{p}}{[i_s]} &+\!\!=& \frac{\partial v_y}{\partial y}{[i_s]}
		\end{array}
		$

		\If{$i_s$ is at source location}
		{
			$
			\arraycolsep=1pt\def\arraystretch{1.75}
			\begin{array}{rrl}		
\nl\hfill	{R^{p}}{[i_s]} &+\!\!=& {s_p}^{(i_t)}
			\end{array}
			$		
		}

		$
		\arraycolsep=1pt\def\arraystretch{1.75}
		\begin{array}{rrl}		
\nl\hfill	{R^p}{[i_s]} &/\!=& {\beta^p}{[i_s]} \\
\nl\hfill	   P{[i_s]} &+\!\!=& {R^p}{[i_s]}
		\end{array}
		$
	}
}
\end{algorithm}
\vspace{1em}
\end{minipage}
\qquad
\begin{minipage}{0.475\linewidth}
\vspace{1em}
\begin{algorithm}[H]
\footnotesize
\SetAlCapHSkip{0ex}
\caption{Simulation procedure with compensated sum.}
\label{alg_compensated_sum}
\For {$i_t = 0 \cdots N_t$}
{
	\For {$i_s = 0 \cdots N^{v_x}_s$}
	{
		$
		\arraycolsep=1pt\def\arraystretch{1.75}
		\begin{array}{rrl}
\nl\hfill	{R^{v_x}}{[i_s]} &\textcolor{kwdcolor}{+\!\!=}& \frac{\partial p}{\partial x}{[i_s]} \\
\nl\hfill	{R^{v_x}}{[i_s]} &/\!=& {\rho^{v_x}}{[i_s]} \\
\nl\hfill	\textcolor{kwdcolor}{\texttt{sum\_3op}} & & \hspace{-1.25em} ( {V_x}{[i_s]} \,, {R^{v_x}}{[i_s]} \,, {V_x}{[i_s]} \,, {R^{v_x}}{[i_s]} )
		\end{array}
		$
	} 

	\For {$i_s = 0 \cdots N^{v_y}_s$}
	{
		$
		\arraycolsep=1pt\def\arraystretch{1.75}
		\begin{array}{rrl}
\nl\hfill	{R^{v_y}}{[i_s]} &\textcolor{kwdcolor}{+\!\!=}& \frac{\partial p}{\partial y}{[i_s]} \\
\nl\hfill	{R^{v_y}}{[i_s]} &/\!=& {\rho^{v_y}}{[i_s]} \\
\nl\hfill	\textcolor{kwdcolor}{\texttt{sum\_3op}} & & \hspace{-1.25em} ( {V_y}{[i_s]} \,, {R^{v_y}}{[i_s]} \,, {V_y}{[i_s]} \,, {R^{v_y}}{[i_s]} )
		\end{array}
		$
	}

	\For {$i_s = 0 \cdots N^{p}_s$}
	{
		$
		\arraycolsep=1pt\def\arraystretch{1.75}
		\begin{array}{rrl}
\nl\hfill	{R^{p}}{[i_s]} &\textcolor{kwdcolor}{+\!\!=}& \frac{\partial v_x}{\partial x}{[i_s]} \\
\nl\hfill	{R^{p}}{[i_s]} &+\!\!=& \frac{\partial v_y}{\partial y}{[i_s]}
		\end{array}
		$

		\If{$i_s$ is at source location}
		{
			$
			\arraycolsep=1pt\def\arraystretch{1.75}
			\begin{array}{rrl}		
\nl\hfill	{R^{p}}{[i_s]} &+\!\!=& {s_p}^{(i_t)}
			\end{array}
			$		
		}

		$
		\arraycolsep=1pt\def\arraystretch{1.75}
		\begin{array}{rrl}		
\nl\hfill	{R^p}{[i_s]} &/\!=& {\beta^p}{[i_s]} \\
\nl\hfill	\textcolor{kwdcolor}{\texttt{sum\_3op}} & & \hspace{-1.25em} ( {P}{[i_s]} \,, {R^{p}}{[i_s]} \,, {P}{[i_s]} \,, {R^{p}}{[i_s]} )
		\end{array}
		$
	}
}
\end{algorithm}
\vspace{1em}
\end{minipage}

After replacing the standard floating point sum at lines 3, 6, and 11 of Algorithm \ref{alg_basic} with the compensated sum (using \texttt{sum\_3op} as an example), we obtain Algorithm \ref{alg_compensated_sum}. 
For example, at line 3 of Algorithm \ref{alg_compensated_sum}, $V_x[i_s]$ stores the floating point sum of $V_x[i_s]$ and $R^{V_x}[i_s]$ at the exit of \texttt{sum\_3op}, while $R^{V_x}[i_s]$ keeps track of the lost bits, which is later applied as compensation at the next time step (line 1).
The other modifications can be interpreted similarly.

\begin{figure}[H]
\centering
\includegraphics[scale=0.225]{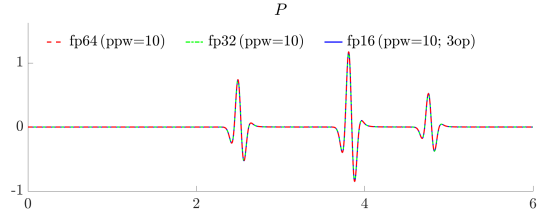}
\captionsetup{width=.9\linewidth,font=small}
\caption{Time history of $P$ at the receiver location.}
\label{fp16_3op_P}
\end{figure}

\begin{figure}[H]
\centering
\includegraphics[scale=0.225]{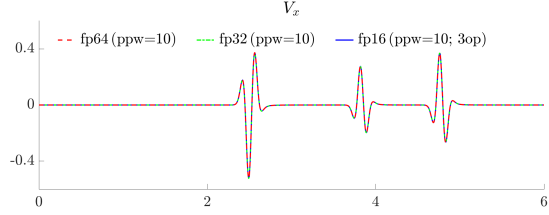}
\captionsetup{width=.9\linewidth,font=small}
\caption{Time history of $V_x$ at the receiver location.}
\label{fp16_3op_Vx}
\end{figure}

\begin{figure}[H]
\centering
\includegraphics[scale=0.225]{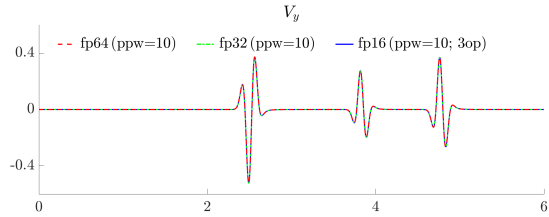}
\captionsetup{width=.9\linewidth,font=small}
\caption{Time history of $V_y$ at the receiver location.}
\label{fp16_3op_Vy}
\end{figure}

\begin{figure}[H]
\centering
\includegraphics[scale=0.225]{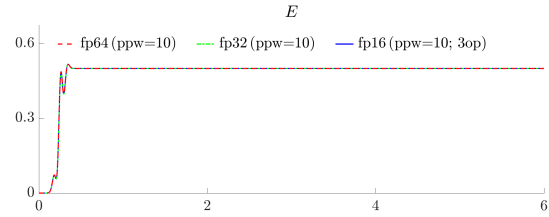}
\captionsetup{width=.9\linewidth,font=small}
\caption{Time history of the overall energy.}
\label{fp16_3op_E}
\end{figure}

\begin{figure}[H]
\centering
\includegraphics[scale=0.225]{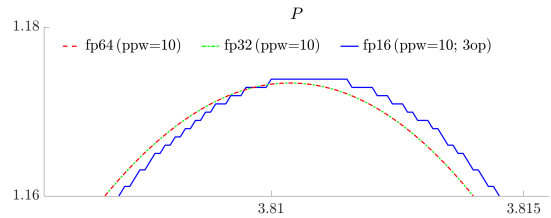}
\captionsetup{width=.9\linewidth,font=small}
\caption{Zoom-in plot of Figure \ref{fp16_3op_P} reveals minute difference in the fp16 simulation result from the fp32 and fp64 results. Comparison with Figure \ref{fp64_fp32_P_zoom_in} indicates that the difference is still less than or on the same level as that resulted from the discretization error.}
\label{fp16_3op_P_zoom_in}
\end{figure}

Below, we present numerical examples to demonstrate that applying compensated sum can drastically improve the quality of the wave simulation. 
The simulation configuration is the same as that used in the previous section.
Simulation results using Algorithm \ref{alg_compensated_sum} are shown in Figures \ref{fp16_3op_P} - \ref{fp16_3op_E} (cf. Figures \ref{fp16_naive_P} - \ref{fp16_naive_E} from the previous section).
From these figures, we observe that the fp16 simulation results closely agree with the fp32 and fp64 results. 
Minute difference can be observed only when zoomed in closely, as shown in Figure \ref{fp16_3op_P_zoom_in} for pressure, which is zoomed in \string~131 times in the vertical direction from Figure \ref{fp16_3op_P}.
Comparing Figure \ref{fp16_3op_P_zoom_in} with Figure \ref{fp64_fp32_P_zoom_in} (which are plotted with the same scale) suggests that this difference is still less than or on the same level as that resulted from the discretization error.
In other words, if the discretization error in the fp64 simulation at ppw=10 is deemed acceptable, the roundoff error accumulated in the fp16 simulation at ppw=10 should be deemed acceptable as well.

\renewcommand{\folderPath}{\string"\string./figure/fp16_3op/\string"}

Before proceeding to the more broad discussions, we make a note here on the storage for the compensation term. 
In the two algorithms above, we use a right hand side term to store the quantity to be added to each solution variable during its update, e.g., $R^{V_x}$ is the right hand side term for $V_x$. 
In Algorithm \ref{alg_compensated_sum}, these right hand side terms serve double duty as the storage for the compensation terms as well, thus there is no need for allocating additional storage.

Careful readers may have noticed that the right hand side terms in Algorithm \ref{alg_basic} could have been eliminated if the stencil application and solution update are combined together.
This is possible because of three special circumstances in this particular case. 
First, the updating equation for each solution variable does not depend on its own spatial derivatives. 
Second, the time integration scheme involves only the immediate previous time step. 
Third, pressure and velocities are staggered in time.

Any violation of the three circumstances listed above, e.g., when using a different time integration scheme such as Runge-Kutta or multi-step methods, or when using a second order formulation of the equation that involves only the pressure, would necessitate the need for temporary storage such as the right hand side terms 
to avoid prematurely overwriting the solution.
Therefore, in these more general cases, we can always expect to find "free" storage for the compensation terms (see \cite{gill1951process} for an example for the Runge-Kutta case).

\section{Discussions}\label{sec_discussions}
While practitioners have had ample experiences with fp64 and fp32 in wave simulations, simulating with fp16 remains a novelty. 
Issues unexpected at fp64 and fp32 may arise at fp16 due to the limited number of bits in the number representation.
In the above, we illustrated a particular issue associated with precision and provided remedy accordingly.
Below, we outline two other potential issues to make the practitioners aware and to possibly spark further investigations from the community.

{\bf Range.} Aside from the reduced precision, the fp16 format is also limited in its range of representable numbers (from $6.5504\times 10^{4}$ to $6.1035\times 10^{-5}$ for normal numbers and to $5.9605\times 10^{-8}$ if subnormals are included; see Table \ref{tbl_formats}). 
Special care may need to be taken to choose the proper units or rescale the problem so that parameters and solution variables remain within the representable range.
To give an example, if using the SI units, the density and wave-speed for water would be $\sim$1000~
$\text{kg}/\text{m}^3$ and $\sim$1500~$\text{m}/\text{s}$, respectively, which would lead to a bulk modulus with value $2.25 \times 10^9$ at unit $\text{kg}/\text{ms}^2$ that exceeds the representable range of the fp16 format.
%

{\bf Operator.} While we focused on the effect of roundoff errors accumulated during the course of the simulation, there is another mechanism roundoff errors can impact the simulated solution, i.e., by inducing errors in the discretization operators (e.g., the stencil). 
The design of discretization operators typically assumes infinite precision.
With limited precision at fp16, the operators actually represented on computer may deviate sufficiently from the intended infinite precision ones to lose approximation property.
This is more likely to become an issue, e.g., when long stencils with a wide range of coefficients or complex time integration schemes are employed.

\section{Conclusion}\label{sec_conclusion}
We have examined the viability of conducting wave simulations using IEEE half precision arithmetic and focused on the impact of the accumulation of roundoff errors on simulation results. 
Degradation in the solution quality when naively switching to half precision are illustrated.
%
The addition operations in the solution update, which form a disguised recursive sum, are identified as the cause of the issue.
Remedy in the form of compensated sum is then provided.
%
We demonstrated that applying compensated sum can greatly restore the solution quality for half precision wave simulations for both the acoustic and elastic cases.
%

Nevertheless, we caution the practitioners that aside from the specific issue addressed here, other issues may still appear in half precision wave simulations, some of which are outlined in the Discussions section. 
These issues, along with errors from discretization, can be convoluted and hard to discern in practice. 
More research work, including carefully designing experiments to expose problems, discerning the causes, and finding remedies, is still needed to achieve a level of confidence in half precision simulations comparable to that in double and single precision simulations.


\section{Acknowledgements}
This research used resources of the Argonne Leadership Computing Facility, a U.S. Department of Energy (DOE) Office of Science user facility at Argonne National Laboratory and is based on research supported by the U.S. DOE Office of Science-Advanced Scientific Computing Research Program, under Contract No. DE-AC02-06CH11357.

\renewcommand*{\theHsection}{\thepart.\thesection}
\appendix

\section{Results with elastic wave equation}\label{App_elastic}
\setcounter{figure}{0}%
%

In this appendix, we demonstrate that the issues associated with the fp16 simulation of the acoustic wave equation illustrated in Section \ref{sec_problem} also appear in the elastic case and that the remedy presented in Section \ref{sec_compensated_sum} still works.

To have some variety, the numerical experiments for the elastic case are tested on NVIDIA A100 GPU (instead of Intel GPU for the acoustic case), where the fp16 type can be accessed through the CUDA type \texttt{\_\_half} defined in header \texttt{cuda\_fp16.h}.

The 2D elastic wave equation under consideration is outlined below, where $v_x$ and $v_y$ are the particle velocities, $\sigma_{xx}$, $\sigma_{xy}$, and $\sigma_{yy}$ are the stress components; $\lambda$ and $\mu$ are Lam\'e parameters that relate with the density $\rho$, compressional wave-speed $c_p$, and shear wave-speed $c_s$ via formulas $\lambda = \rho (c_p^2 - 2 c_s^2)$ and $\mu = \rho c_s^2$.

\begin{equation}
\arraycolsep=2.5pt\def\arraystretch{2.25}
\label{elastic_wave_equation}
\left\{
\arraycolsep=2.5pt\def\arraystretch{1.5}
\begin{array}{rcl}
\displaystyle \frac{\partial v_x}{\partial t} &=& \displaystyle \frac{1}{\rho} \left( \frac{\partial \sigma_{xx}}{\partial x} + \frac{\partial \sigma_{xy}}{\partial y} \right)
\\[2ex]
\displaystyle \frac{\partial v_y}{\partial t} &=& \displaystyle \frac{1}{\rho} \left( \frac{\partial \sigma_{xy}}{\partial x} + \frac{\partial \sigma_{yy}}{\partial y} \right) + s^{v_y}
\\[2ex]
\displaystyle \frac{\partial \sigma_{xx}}{\partial t} &=& \displaystyle \left(\lambda + 2\mu\right) \frac{\partial v_x}{\partial x} + \lambda \frac{\partial v_y}{\partial y}
\\[2ex]
\displaystyle \frac{\partial \sigma_{xy}}{\partial t} &=& \displaystyle \mu \frac{\partial v_y}{\partial x} + \mu \frac{\partial v_x}{\partial y}
\\[2ex]
\displaystyle \frac{\partial \sigma_{yy}}{\partial t} &=& \displaystyle \lambda \frac{\partial v_x}{\partial x} + \left(\lambda + 2\mu\right) \frac{\partial v_y}{\partial y}
\end{array}
\right.
\end{equation}

\begin{figure}[H]
\centering
\includegraphics[scale=0.2]{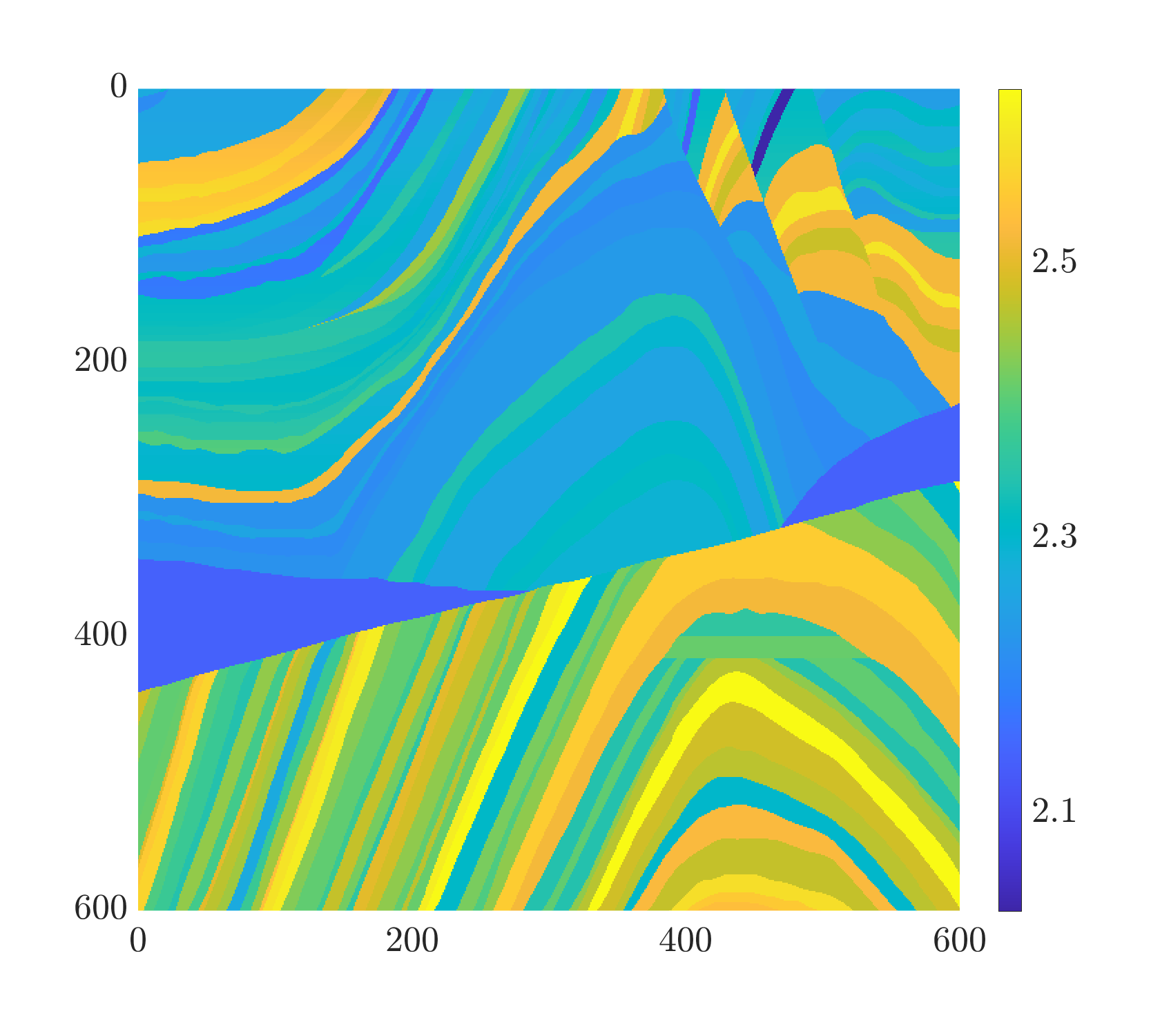}
\vspace{-1.5em}
\caption{Illustration of the medium parameter (compressional wave-speed) used in the elastic simulations.}
\label{medium_parameter_vp}
\end{figure}

The simulation configuration is mostly the same as the acoustic case, with a few exceptions outlined below.
A singular source $s^{v_y}$ on the vertical velocity $v_y$ is considered here, which models a source generated by, e.g., a weight dropping mechanism. 
Periodic boundary condition is considered for the $x$-direction and free surface boundary condition is considered for the $y$-direction here.
The medium parameters are taken from a subsection of the Marmousi2 model (\cite{martin2006marmousi2}) with the minimal shear wave-speed at $\sim$1.0117 km/s, the maximal compressional wave-speed at $\sim$4.6992 km/s, and density ranges from $\sim$2.0293 to $\sim$2.6230 $\text{Gt}/\text{km}^3$. The compressional wave-speed is illustrated in Figure \ref{medium_parameter_vp}.
Both source and receiver are placed at 10 grid spacings below the top surface in order to generate and record the surface wave, at 200th and 400th grid point horizontally, respectively.

The simulation results without compensation are illustrated in Figures \ref{fp16_elastic_Vy_0op} - \ref{fp16_elastic_E_0op} for $V_y$ at the receiver location, its zoom-in plot for the last second, and the overall energy evolution. 
(The other solution components are omitted here to conserve space.)
We can see from these figures that the same issues of small wiggles and energy loss also appear here. 
In comparison, the simulation results with compensation are illustrated in Figures \ref{fp16_elastic_Vy_6op} - \ref{fp16_elastic_E_6op}. 
We can readily observe from these figures that by applying the compensated sum, the simulation quality is markedly improved.
(We used the \texttt{sum\_6op} version here to have some variety. The \texttt{sum\_3op} gives very similar results.)

\renewcommand{\folderPath}{\string"\string./figure/fp16_6op_elastic/\string"}

\begin{figure}[H]
\centering
\includegraphics[scale=0.225]{\string"\folderPath/Vx_0op\string".png}
\captionsetup{width=.9\linewidth,font=small}
\caption{Time history of $V_y$ at the receiver location.}
\label{fp16_elastic_Vy_0op}
\end{figure}

\begin{figure}[H]
\centering
\includegraphics[scale=0.225]{\string"\folderPath/Vx_zoom_in_0op\string".png}
\captionsetup{width=.9\linewidth,font=small}
\caption{Zoom-in plot of Figure \ref{fp16_elastic_Vy_0op}.}
\label{fp16_elastic_Vy_0op_zoom_in}
\end{figure}

\begin{figure}[H]
\centering
\includegraphics[scale=0.225]{\string"\folderPath/E_0op\string".png}
\captionsetup{width=.9\linewidth,font=small}
\caption{Time history of the overall energy.}
\label{fp16_elastic_E_0op}
\end{figure}

\begin{figure}[H]
\centering
\includegraphics[scale=0.225]{\string"\folderPath/Vx_6op\string".png}
\captionsetup{width=.9\linewidth,font=small}
\caption{Time history of $V_y$ at the receiver location.}
\label{fp16_elastic_Vy_6op}
\end{figure}

\begin{figure}[H]
\centering
\includegraphics[scale=0.225]{\string"\folderPath/Vx_zoom_in_6op\string".png}
\captionsetup{width=.9\linewidth,font=small}
\caption{Zoom-in plot of Figure \ref{fp16_elastic_Vy_6op}.}
\label{fp16_elastic_Vy_6op_zoom_in}
\end{figure}

\begin{figure}[H]
\centering
\includegraphics[scale=0.225]{\string"\folderPath/E_6op\string".png}
\captionsetup{width=.9\linewidth,font=small}
\caption{Time history of the overall energy.}
\label{fp16_elastic_E_6op}
\end{figure}

\section{Experiments on stencil application}\label{App_promotion}
\setcounter{figure}{0}%

In the following, we demonstrate that stencil application, despite involving majority of the floating point operations in Algorithm \ref{alg_basic}, is not responsible for the undesirable features in the fp16 simulation results illustrated in Figures \ref{fp16_naive_P} - \ref{fp16_naive_E}.

In Figures \ref{fp16_naive_promotion_P} and \ref{fp16_naive_promotion_E}, we show the simulation results using Algorithm \ref{alg_basic}, but with the stencil application promoted to fp64 and the solution update retained at fp16. 
As can be readily seen, the undesirable wiggles and energy loss are still present in the simulation results.
On the other hand, once we include compensation for the solution update, i.e., switch to Algorithm \ref{alg_compensated_sum}, these issues are addressed as shown in Figures \ref{fp16_3op_promotion_P} and \ref{fp16_3op_promotion_E}, validating the claim that lack of precision at the solution update is the cause underlying these issues.
(To conserve space, the plots for $V_x$ and $V_y$ are omitted here.)

\renewcommand{\folderPath}{\string"\string./figure/fp16_naive_promotion/\string"}

\begin{figure}[H]
\centering
\includegraphics[scale=0.225]{\string"\folderPath/P\string".png}
\captionsetup{width=.9\linewidth,font=small}
\caption{Time history of $P$ at the receiver location.}
\label{fp16_naive_promotion_P}
\end{figure}

\begin{figure}[H]
\centering
\includegraphics[scale=0.225]{\string"\folderPath/E\string".png}
\captionsetup{width=.9\linewidth,font=small}
\caption{Time history of the overall energy.}
\label{fp16_naive_promotion_E}
\end{figure}

\renewcommand{\folderPath}{\string"\string./figure/fp16_3op_promotion/\string"}

\begin{figure}[H]
\centering
\includegraphics[scale=0.225]{\string"\folderPath/P\string".png}
\captionsetup{width=.9\linewidth,font=small}
\caption{Time history of $P$ at the receiver location.}
\label{fp16_3op_promotion_P}
\end{figure}

\begin{figure}[H]
\centering
\includegraphics[scale=0.225]{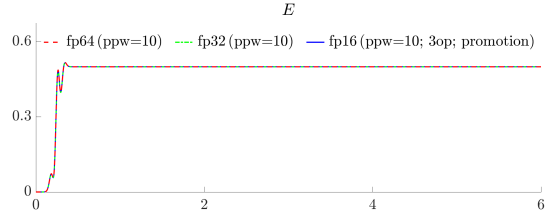}
\captionsetup{width=.9\linewidth,font=small}
\caption{Time history of the overall energy.}
\label{fp16_3op_promotion_E}
\end{figure}

The fact that roundoff errors occurred at different operations have different impact on the final outcome is not surprising. 
To give a pedagogical example, consider the following procedure to calculate the inner product of two vectors, i.e., $s =x^T y$, where $x_i$ and $y_i$ are the components of $x$ and $y$, respectively, and $s_i$ is the partial outcome at step $i$, with $s_n$ being the final outcome that approximates $s$.

\begin{equation}
\label{inner_product_procedure}
\arraycolsep=2pt
\begin{array}{lcrll}
s_0 & = & & & (x_0 \otimes y_0) \\
s_1 & = & s_0 & \oplus & (x_1 \otimes y_1) \\
s_2 & = & s_1 & \oplus & (x_2 \otimes y_2) \\
& \vdots & \\
s_n & = & s_{n-1} & \oplus & (x_{n} \otimes y_{n})
\end{array}
\end{equation}

In the above procedure, $\otimes$ and $\oplus$ denote the floating point multiplication and addition, respectively.
There are roughly the same amount of $\otimes$ and $\oplus$ in the above procedure.
We assume that they satisfy the following relations with the exact multiplication and addition ($u$ stands for unit roundoff), which is guaranteed if the IEEE standard is followed.
\begin{equation}
\label{inner_product_relations}
\arraycolsep=2pt
\begin{array}{ccc}
(x_i \otimes y_i) &=& (x_i \times y_i)(1 + \epsilon_i)\,, \quad \vert\epsilon_i\vert \leq u \\
(s_i \oplus  y_i) &=& (s_i +  	  y_i)(1 + \delta_i  )\,, \quad \vert\delta_i\vert   \leq u \\
\end{array}
\end{equation}
Substituting these relations into the procedure above, we obtain the following expression for $s_n$.
\begin{equation}
\label{inner_product_outcome}
\arraycolsep=2pt
\begin{array}{lcl}
s_n & = & (x_0 \times y_0)(1 + \textcolor{black}{\epsilon}_0)\prod_{i=1}^n(1 + \textcolor{black}{\delta}_i) \\[0.5ex]
    & + & (x_1 \times y_1)(1 + \textcolor{black}{\epsilon}_1)\prod_{i=1}^n(1 + \textcolor{black}{\delta}_i) \\[0.5ex]
    & + & (x_2 \times y_2)(1 + \textcolor{black}{\epsilon}_2)\prod_{i=2}^n(1 + \textcolor{black}{\delta}_i) \\
    & \vdots & \\
    & + & (x_n \times y_n)(1 + \textcolor{black}{\epsilon}_n)\prod_{i=n}^n(1 + \textcolor{black}{\delta}_i) \\ 
\end{array}
\end{equation}

From equation \ref{inner_product_outcome}, we see that the roundoff errors from addition ($\delta_i$) compound together while those from multiplication ($\epsilon_i$) do not.
For example, the product $x_0 \times y_0$ is perturbed by all $\delta_i$ committed at the subsequent additions, but only by $\epsilon_0$ committed at this multiplication itself.
In this case, the roundoff errors committed at $\oplus$ are problematic as they compound, 
while those committed at $\otimes$ are benign.


%

\bibliographystyle{alpha}
\bibliography{refs}

\end{document}